\documentclass{amsart}
\usepackage{tipa}
\usepackage{amssymb}
\usepackage{stmaryrd}
\usepackage{graphicx}

\newtheorem{theorem}{Theorem}[section]

\newtheorem{problem}[theorem]{Problem}

\newtheorem{proposition}[theorem]{Proposition}
\newtheorem{corollary}[theorem]{Corollary}

\theoremstyle{definition}

\theoremstyle{remark}
\newtheorem{remark}[theorem]{Remark}

\numberwithin{equation}{section}

\begin{document}

\author{TIANXIN CAI}
\address{School of Mathematical Sciences, Zhejiang University, Hangzhou 310027, People's Republic of China }

\email{txcai$@$zju.edu.cn}

\thanks{This research was supported by the National Natural Science Foundation of China (Grant No.~11571303).}

\author{YONG ZHANG}
\address{School of Mathematical Sciences, Zhejiang University, Hangzhou 310027, People's Republic of China }

 \email{zhangyongzju$@$163.com}

\title[Congruent numbers on the right trapezoid]{Congruent numbers on the right trapezoid}



\begin{abstract}
We introduce and study a new kind of congruent number problem on the
right trapezoid.
\end{abstract}
\maketitle

\section{Introduction}
A congruent number is a positive integer that is the area of a right
triangle with three rational number sides. In other words, $n$ is a
congruent number if and only if there is a right triangle with
rational sides $a,b,c\in\mathbb{Q}^+$ such that
\[a^2+b^2=c^2, ab=2n.\]

Let $x=n(a+c)/b,y=2n^2(a+c)/b^2$, we get a family of elliptic curves
\[E_n: y^2=x^3-n^2x.\]
Conversely, if we have rational solutions $(x,y)$ on the elliptic
curve $E_n$ with $y \neq 0$, then the rational numbers
\[a=\bigg|\frac{2nx}{y}\bigg|,b=\bigg|\frac{x^2-n^2}{y}\bigg|,c=\bigg|\frac{x^2+n^2}{y}\bigg|\]
are the sides of a right triangle with area $n$.

Due to the homogeneity of the condition $a^2+b^2=c^2,$ we only need
consider the square-free positive integers. Determining whether a
given square-free positive integer is a congruent is the congruent
number problem. It has not yet been solved in general. Many
mathematicians studied this problem, such as Fibonacci, Fermat and
Euler. For more information about this problem, we can refer to
\cite{Alter-Curtz,Guy,Tian1,Tian2,Top-Yui}.

In particular, there is a generalized congruent number problem. A
positive integer $n$ is called $t$-congruent number (see
\cite{Top-Yui}) if there are positive rational numbers $a,b,c$ such
that
\[a^2=b^2+c^2-2bc\frac{t^2-1}{t^2+1},bc\frac{2t}{t^2+1}=2n.\]
The case $t=1$ corresponds to the classical congruent number
problem. Its corresponding elliptic curve is
\[E_{n,t}: y^2=x(x-\frac{n}{t})(x+nt).\]

Now we consider a new kind of congruent number problem which related
with the area of right trapezoid.

\noindent{\bf Definition 1.1.}\emph{ A positive integer $n$ is
called $i$-congruent number, if it is the area of the right
trapezoid $($Figure $1)$ with $a,b,c \in\mathbb{Z}^+,d\in\mathbb{N},
(b,c)=1.$}

When $d=0$, $n$ is the classic congruent number with $a,b,c \in\mathbb{Z}^+$. By this definition, we have
\begin{equation}n=(a+d)b/2,(a-d)^2+b^2=c^2, (b,c)=1.\end{equation}

\begin{figure}[htbp]
\centering
\includegraphics[scale=0.6]{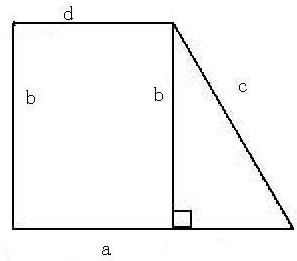}
\caption{}\label{fig:side:a}
\end{figure}

From the property of Pythagoras triples, we have $(a-d,b)=1,$ and
$a-d,b$ have opposite parity. Then
\[(a-d,b)=(2xy,x^2-y^2),~or~(x^2-y^2,2xy),\]
where $x>y,(x,y)=1,$ and $x,y$ have opposite parity. Observing this, we have

\begin{proposition}
A positive integer $n$ is $i$-congruent number if and only if
$n=pk$, where $p$ is any odd prime and $k\geq \frac{p^2-1}{4}$, or
$n=2^ik,k\geq 2^{2i}-1,i\geq 1$, where $k$ is any odd integer.
\end{proposition}

The $i$-congruent numbers in Proposition 1.1 may coincide in some
cases, the intersection set of them less than 100 is
\[\{6,18,30,42,50,54,60,70,78,84,90,98,100\}.\]

On the other hand, $n>1$ is non-$i$-congruent number if and only if
$n$ has the following form
\[\begin{split}p,~p^2(p\neq3),~pq(5<p<q<\frac{p^2-1}{4}),\\
~2^i(i\geq 0),~2^ip(i\geq 2,2^{1+i/2}<p<2^{2i}-1),\end{split}~~~~~~~~(*)\]where $p,q$ are primes.

For examples, the number of non-$i$-congruent numbers less than 100
is 46, where the non-primes have 21 which are
\[1,4,8,16,20,25,28,32,49,52,56,58,62,64,74,77,82,86,88,91,94.\]But we
have

\begin{proposition}Almost every positive integer is $i$-congruent number.\end{proposition}

\begin{theorem}Let $f(x)$ be the number of non-$i$-congruent numbers less than $x$, then
\[f(x)\sim~\frac{cx}{\log x},\] where $c=1+\ln2$.
\end{theorem}

For some $i$-congruent numbers $m$, there might be many right
trapezoids having the area $n$. By Proposition 1.1 and its proof,
for any positive integer $m$, there are infinitely many congruent
integers $n$ such that we have $m$ or more right trapezoids which
have the area $n$, it is only to need take $n=p_1\cdots p_mn',$
where $p_i$ is the $i$-$th$ prime and $n'\geq \frac{(p_m-1)^2}{4}$
is an arbitrary integer.

We find 16 $i$-congruent numbers with $d=0$ less than 1000, i.e.,
\[6,30,60,60,84,96,180,210,330,480,486,504,546,630,840,924,960,\]
where there are two right triangles with
$(a,b,c)=(21,20,29),(35,12,37)$ have the same area $n=210.$ Then we
have

\begin{proposition}Let $g(x)$ be the number of $i$-congruent numbers less than $x$ for $d=0$
with $a>b$ and $a,b \in \mathbb{Z}^+,$ including the repeated cases,
then \[\frac{\sqrt{x}}{2}+O(1)<g(x)\leq \frac{1}{2
\sqrt[3]{4}}x^{\frac{2}{3}}+O(x^{\frac{5}{9}}).\]
\end{proposition}

Next, we give the definitions of two new kinds of congruent numbers.

\noindent{\bf Definition 1.2.}\emph{ A positive integer $n$ is
called $k$-congruent number, if it is the area of the right
trapezoid $($Figure $1)$ with $a,b,c,d \in\mathbb{Q}^+,
k\in\mathbb{Z}^+,k\geq1,a=kd.$}

By this definition, we have
\begin{equation}n=(a+d)b/2,(a-d)^2+b^2=c^2.\end{equation}

\noindent{\bf Definition 1.3.}\emph{ A positive integer $n$ is
called $d$-congruent number, if it is the area of the right
trapezoid $($Figure $2)$ with $a,b,c\in\mathbb{Q}^+,
d\in\mathbb{N}.$}

When $d=0$, $n$ is the classic congruent number. By this definition, we have
\begin{equation}n=(a+2d)b/2,a^2+b^2=c^2.\end{equation}

\begin{figure}[htbp]
\centering
\includegraphics[scale=0.6]{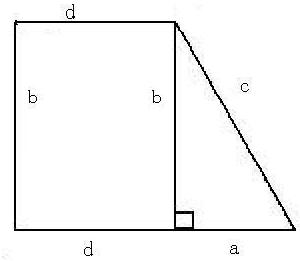}
\caption{}\label{}
\end{figure}

For $k$-congruent number and $d$-congruent number, we have

\begin{proposition} For every positive integer $n$, there exists a $k$ such that $n$ is $k$-congruent number.\end{proposition}

\begin{theorem} For every positive integer $n$, there exists a $d$ such that $n$ is $d$-congruent number.\end{theorem}

In section 2, we will give the proofs of propositions and theorems.
We consider the family of elliptic curves about $k$-congruent
numbers in section 3, and prove two propositions about $d$-congruent
numbers in section 4.

\vskip 10pt
\section{Proofs of Propositions and Theorems}
\emph{Proof of Proposition 1.1.} By Definition 1.1, $n$ is an
$i$-congruent number if and only if the Diophantine system
\[a+d=\frac{2n}{x^2-y^2},a-d=2xy,~or~a+d=x^2-y^2,a-d=\frac{n}{xy}\]
has non-negative integer solutions.

Solving them, we get
\[a=\frac{n}{x^2-y^2}+xy,d=\frac{n}{x^2-y^2}-xy,~or~2a=\frac{n}{xy}+x^2-y^2,2d=\frac{n}{xy}-x^2-y^2.\]
Hence, when $a-d$ be even, $a,d$ are non-negative integers if and
only if $n=(x^2-y^2)k,k\geq xy,$ and when $a-d$ be odd $a,d$ are
non-negative integers if and only if $n=(x^2-y^2)k,k\geq xy.$

Because every odd prime can be presented by the difference of two
coprime integers, so the positive integer $n=pk$ is $i$-congruent
number, where $p$ is odd prime and $k\geq \frac{p^2-1}{4}.$ Next, we
will prove the opposite direction is also right. Let
$n=(x^2-y^2)k,k\geq xy,x>y,(x,y)=1$, and $x,y$ have opposite parity.
If $p|x^2-y^2,$ put $x^2-y^2=ps,$ then $n=p(sk)$. If $p>s$, let
\[x=\frac{p+s}{2},y=\frac{p-s}{2},xy=\frac{p^2-s^2}{4},\]hence
$sk\geq s\frac{p^2-s^2}{4}\geq \frac{p^2-1}{4}.$ If $p<s$, let
 \[x=\frac{p+s}{2},y=\frac{s-p}{2},xy=\frac{s^2-p^2}{4},\]
 hence $sk\geq s\frac{s^2-p^2}{4}\geq \frac{p^2-1}{4}.$ If $p=s$, let
 \[x=\frac{p^2+1}{2},y=\frac{p^2-1}{2},xy=\frac{p^4-1}{4},\]hence $sk\geq p\frac{p^4-1}{4}\geq \frac{p^2-1}{4}.$
So in any case, we have $sk \geq \frac{p^2-1}{4}.$ Let $k'=ks,$ we
get $n=pk',k'\geq \frac{p^2-1}{4}.$

By the same method, the positive integer $n=2^ik,k\geq
2^{2i}-1,i\geq 1$ is $i$-congruent number, where $k$ is odd integer,
only if we let $x=2^i,y=1$ in $n=xyk.$ Next, we want to prove the
opposite direction is also right. Let $n=xyk,k\geq
x^2-y^2,x>y,(x,y)=1$, and $x,y$ have opposite parity. If $2^i||xy,$
we can put $xy=2^it,2\nmid t,$ then $n=2^itk$. If $2^i>t$, let
\[x=2^i,y=t,x^2-y^2=2^{2i}-t^2,\]
hence $tk\geq t(2^{2i}-t^2)\geq 2^{2i}-1.$ If $2^i<t$, let
\[x=t,y=2^i,x^2-y^2=t^2-2^{2i},\]hence
$tk\geq t(t^2-2^{2i})\geq 2^{2i}-1.$ So in any case, we have $tk
\geq 2^{2i}-1.$ Let $k'=kt,$ we get $n=2^ik',k'\geq 2^{2i}-1,$ where $k'$ is an odd integer. \hfill $\Box$ \\

\emph{Proof of Proposition 1.2.} We only need to prove the case for
\[a=\frac{n}{x^2-y^2}+xy,d=\frac{n}{x^2-y^2}-xy.\]It's
easy to see that for a fix odd prime, when $x$ is large enough, the
number of $i$-congruent numbers $n$ in the interval $[1,x]$
satisfying
\[1\leq n=pk\leq x, k\geq \frac{p^2-1}{4}\]
has the proportion $1/p.$ By the Pigeonhole principle, taking $p$ be
the first $s$ primes, under the above condition, the number of
$i$-congruent numbers in the interval $[1,x]$ has the proportion
 \[1-\prod_{j=1}^{s}(1-\frac{1}{p_j}),\]where $p_j$ is the $j$-$th$ prime.
Since \[\lim_{s\rightarrow
\infty}(1-\prod_{j=1}^{s}(1-\frac{1}{p_j}))=1,\]
 hence almost every positive integer is $i$-congruent number. \hfill $\Box$ \\

\emph{Proof of Theorem 1.3.} It's easy to see that the number of the
$2$-$th$, $4$-$th$ and $5$-$th$ term in $(*)$ are $O(\sqrt{x})$,
$O(\log x)$ and $O(x^{\frac{2}{3}})$. Let $\pi(x)$ be the prime
function, the number of the first term in $(*)$ is
\[\pi(x)=\frac{x}{\log x}+O\bigg(\frac{x}{\log^2x}\bigg).\]

In the following, we consider the third term, i.e.,
\[\sum_{pq\leq x,5<p<q<p^2/4}1=\sum_{5<p\leq \sqrt[3]{4x}}~~\sum_{p<q<p^2/4}1+\sum_{\sqrt[3]{4x}<p\leq \sqrt{x}}~\sum_{p<q<x/p}1=\sum_1+\sum_2.\]
In view of $\frac{x^2}{\log x}$ is an increasing function, then
\[\sum_1=O\bigg(\sum_{5<p\leq \sqrt[3]{4x}}(\pi(p^2/4)-\pi(p))\bigg)=O\bigg(\sum_{p\leq \sqrt[3]{4x}}\frac{p^2}{\log p}\bigg)=O\bigg(\frac{x}{\log^2 x}\bigg).\]
And\[\begin{split}\sum_2&=\sum_{\sqrt[3]{4x}< p\leq \sqrt{x}}(\pi(x/p)-\pi(p))\\
&=\sum_{\sqrt[3]{x}< p\leq \sqrt{x}}\pi(x/p)-\sum_{\sqrt[3]{x}< p\leq \sqrt[3]{4x}}\pi(x/p)-\sum_{\sqrt[3]{4x}< p\leq \sqrt{x}}\pi(p)\\
&=\sum_3-\sum_4-\sum_5.\end{split}\]
Similarly, by the prime number theorem and the identity
\[\sum_{p\leq x}\frac{1}{p}=\log\log x+c_1+O\bigg(\frac{1}{\log x}\bigg),\]where $c_1$ is a constant, we get
\[\sum_4=\sum_5=O\bigg(\frac{x}{\log^2 x}\bigg).\]
At last, we estimate $\sum\limits_{3}.$ Let
\[c_n=
\begin{cases}
1,~n~is~a~pime, \\
0, ~n~is~not~a~pime,
\end{cases}\]
and \[f(n)=\frac{x/n}{\log(x/n)}.\]By the prime number theorem and
the sum formula of Abel, we get
\[\begin{split}&\sum_{\sqrt[3]{x}< p\leq \sqrt{x}}\frac{x/p}{\log(x/p)}=\sum_{\sqrt[3]{x}< n\leq \sqrt{x}}\frac{x/n}{\log (x/n)}c_n\\
&=\pi(\sqrt{x})f(\sqrt{x})-\pi(\sqrt[3]{x})f(\sqrt[3]{x})- \int_{\sqrt[3]{x}}^{\sqrt{x}}\pi(x)f'(t)dt \\
&=\int_{\sqrt[3]{x}}^{\sqrt{x}}\frac{x}{t\log t\log
\frac{x}{t}}dt+O\bigg(\frac{x}{\log^2 x}\bigg).
\end{split}\]
Taking $t=x^{u}$ in the above integral, we have
\[\begin{split} \sum_3&=\frac{x}{\log x}\int_{\frac{1}{3}}^{\frac{1}{2}}\frac{du}{u(1-u)}+O\bigg(\frac{x}{\log^2x}\bigg)\\
&=\frac{x\ln2}{\log x}+O\bigg(\frac{x}{\log^2 x}\bigg).\end{split}\]
Combining the above sum formulas, the number of non-$i$-congruent
numbers less than $x$ is
\[\frac{cx}{\log x},\] where $c=1+\ln2$. \hfill $\Box$ \\

\emph{Proof of Proposition 1.4.} It's easy to see that for $d=0$ the
$i$-congruent numbers have the form $st(s^2-t^2)$, where $s>t>0,$
$\frac{s}{(s,t)}$ and $\frac{t}{(s,t)}$ has the opposite parity. Let
us consider
\[st(s^2-t^2)\leq x,\]because of the repeated cases and the parity, the number of the pair $(s,t)$ satisfying the above inequity great than $g(x)$.
Noting that \[st(s^2-t^2)\geq(t+1)t(2t+1)>2t^3,\] then
$t<\sqrt[3]{x/2}$. From\[st(s^2-t^2)>t(s-t)^3,\] we have $s\leq
\sqrt[3]{x/t}+t$. Then
\[g(x)\leq \sum_{t=1}^{\sqrt[3]{x/2}}(\sqrt[3]{x/t}+t)\leq \frac{1}{2 \sqrt[3]{4}}x^{\frac{2}{3}}+O(x^{\frac{5}{9}}).\]
On the other hand, for any integer $t$, taking $s=[\sqrt[3]{x/t}],$
we have \[st(s^2-t^2)<s^3t\leq x.\]
Then
\[g(x)\leq \frac{1}{2}\sum_{1\leq t\leq \sqrt[4]{x}}([\sqrt[3]{x/t}]-t)=\frac{\sqrt[2]{x}}{2}+O(1).\]
Combining these, we complete the proof of Proposition 1.4. \hfill $\Box$ \\

Next, we use some facts about congruent numbers to prove Proposition
1.5 and the theory of elliptic curve to prove Theorem 1.6.\\

\emph{Proof of Proposition 1.5.} From (1.2), we know that for $k=1$,
the right trapezoid degenerates into a rectangle. It's easy to see
that there are \[a=d=\frac{2}{t},c=b=\frac{t}{2}\] such that the
area is 1, then 1 is a $k$-congruent number. This is a trivial case.
In the following, we consider the case $k\geq 2$.

Let \[b=\bigg|\frac{x^2-(k^2-1)^2n^2}{(k+1)y}\bigg|,d=\bigg|\frac{2nx}{y}\bigg|,\] we get a family of
elliptic curves
\[E_{n,k}: y^2=x^3-(k^2-1)^2n^2x,\]and
\[a=\bigg|\frac{2knx}{y}\bigg|,c=\bigg|\frac{x^2+(k^2-1)^2n^2}{(k+1)y}\bigg|.\]
For a given $k\geq 2$, $E_{n,k}$ is the special case of congruent
number curve, we call it $k$-congruent number curve. Noting that
$n^3-n$ is a congruent number, let $k=n$, then $E_{n,n}$ has
positive rank, which leads to $n$ is a $n$-congruent number.

Therefore, for every positive integer $n$, there exists a $k$ such that $n$ is $k$-congruent number. \hfill $\Box$\\

For examples, when $k=n=2$, from $E_{2,2}$, we have
\[(a,b,c,d)=\bigg(\frac{8}{3},1,\frac{5}{3},\frac{4}{3}\bigg),\bigg(\frac{80}{7},\frac{7}{30},\frac{1201}{210},\frac{40}{7}\bigg),
\bigg(\frac{6808}{4653},\frac{1551}{851},\frac{7776485}{3959703},\frac{3404}{4653}\bigg)\]
such that 2 is a $2$-congruent number.

When $k=n=3$, from $E_{3,3}$, we have \[(a,b,c,d)=\bigg(\frac{9}{4},2,\frac{5}{2},\frac{3}{4}\bigg),\bigg(\frac{21}{40},\frac{60}{7},\frac{1201}{140},\frac{7}{40}\bigg),
\bigg(\frac{851}{517},\frac{4653}{1702},\frac{7776485}{2639802},\frac{851}{1551}\bigg)\] such that 3 is a $3$-congruent number.\\

\emph{Proof of Theorem 1.6.} Let
 \[\begin{cases}\begin{split}
&a=\frac{(3x-d^2-3n)(3x-d^2+3n)}{3(-3y+3dx-d^3)},\\
&b=\frac{2n(3x-d^2)}{(-3y+3dx-d^3)},\\
&c=\frac{(9-6d^2)x^2+9n^2+d^4}{3(-3y+3dx-d^3)}
\end{split}\end{cases}\]
in (1.3), we get a family of elliptic curves \[E_{n,d}:
y^2=x^3-\frac{3n^2+d^4}{3}x+\frac{(9n^2+2d^4)d^2}{27}.\] We call it
$d$-congruent number curve. To prove for every positive integer $n$
there are $a,b,c,d$ satisfying (1.3), let $d=3n,$ we get
\[E_{n,3n}: y^2=x^3-(1+27n^2)n^2x+3n^4(1+18n^2).\]
The discriminant of $E_{n,3n}$ is $\Delta=(4+81n^2)n^6$. When
$n\geq1,$ we have $\Delta>0,$ this means that $E_{n,3n}$ is
nonsingular.

To find a solution $a,b,c,d$ satisfying (1.3), we need to find a
suitable point on $E_{n,3n}$. Noting that the point $P=(-6n^2,3n^2)$
lies on $E_{n,3n}$, using the group law on the elliptic curve, we
obtain the point
\[[2]P=\bigg(\frac{(27n^2+1)(243n^2+1)}{36},-\frac{(81n^2+1)(6561n^4+324n^2-1)}{216}\bigg).\]
It's easy to see that the $x$-coordinate of the point $[2]P$ is not
in $\mathbb{Z}$ for every $n\geq1$. By the Nagell-Lutz Theorem (see
p. 56 of \cite{Silverman-Tate}), for all $n\geq4$ the point $[2]P$
is of infinite order. Then there are infinitely many rational points
on $E_{n,3n}$. Moreover, the point $[2]P$ is the exact point such
that $a,b,c,d=3n$ satisfying (1.3), which leads to
 \[\begin{cases}\begin{split}
&a=\frac{(729n^3-81n^2+27n+1)(9n-1)}{6(1+81n^2)},\\
&b=\frac{12n(1+81n^2)}{(1+9n)(729n^3+81n^2+27n-1)},\\
&c=\frac{43046721n^8+2125764n^6+39366n^4+1620n^2+1}{6(1+81n^2)(1+9n)(729n^3+81n^2+27n-1)}.
\end{split}\end{cases}\]
For every $n\geq1$, we have $a,b,c\in \mathbb{Q}^+$.

Therefore, for every positive integer $n$, there exists a $d$ such that $n$ is $d$-congruent number. \hfill $\Box$\\

For examples, when $n=1,2,3$, from the elliptic curve $E_{n,3n}$, we
have
\[(a,b,c,d)=\bigg(\frac{1352}{123},\frac{123}{1045},\frac{1412921}{128535},3\bigg)\]
such that 1 is a $3$-congruent number,
\[(a,b,c,d)=\bigg(\frac{94571}{1950},\frac{7800}{117971},\frac{11156645809}{230043450},6\bigg)\]
such that 2 is a $6$-congruent number and
\[(a,b,c,d)=\bigg(\frac{123734}{1095},\frac{3285}{71722},\frac{8874450677}{78535590},9\bigg)\]
such that 3 is a $9$-congruent number.

\vskip 10pt
\section{Further Consideration on $E_{n,k}$}
In the proof of Proposition 1.5, we get a family of elliptic curves
\begin{equation}
E_{n,k}: y^2=x^3-(k^2-1)^2n^2x.
\end{equation}
$E_{n,k}$ has four integer points \[(x,\pm
y)=(0,0),((k^2-1)n,0),(-(k^2-1)n,0)\] and the point at infinity. In
[1], the author listed some classes of congruent numbers. For
examples,
\begin{equation}
 n=\alpha^4+4\beta^4,2\alpha^4+2\beta^4,\alpha^4-\beta^4
\end{equation}
are congruent numbers for $\alpha,\beta \in \mathbb{Z}^+.$ Then we
have

\begin{proposition} For a fixed $k\geq2$, there are infinitely many positive integers $n$ which are $k$-congruent numbers.
\end{proposition}
\emph{Proof.} Let $(k^2-1)n=\alpha^4-\beta^4$. Put
\[\alpha=k^2,\beta=1,\]
then
\[n=k^2+1.\]
Therefore, for a fixed $k\geq2$, there are infinitely many positive integers $n=k^2+1$ which are $k$-congruent numbers. \hfill $\Box$\\

Next, we consider the problem for a fixed $n$ whether there are
infinitely many $k$ such that $n$ is a $k$-congruent number. Noting
that $m^3-m$ is a congruent number, we consider the Diophantine
equation
\[n(k^2-1)=m^3-m.\]For every $n>1$, this equation has three integer solutions
\[(k,m)=(n,n),(8n-3,4n-1),(8n+3,4n+1).\]Hence, for each $n\geq 2,$ there are at least three $k$ such that $n$ is a $k$-congruent number.

In view of $(k^4-1)u^2$ is a congruent number for each $k>1$, then
there are infinitely many $k$ such that $u^2$ is a $k$-congruent
number. By $(k^2-1)n=\alpha^4-\beta^4$, we consider the Diophantine
systems
\[\lambda(k+1)=\alpha^2-\beta^2,n(k-1)=\lambda(\alpha^2+\beta^2),\]
which lead to the Pell's equation
\[(n-\lambda^2)\alpha^2-(n+\lambda^2)\beta^2=2n\lambda.\]
When $\lambda=1$, it's easy to prove they have infinitely many
integer solutions for $n=2,5$. When $\lambda=2$, the same result
holds for $n=10,13,52$. When $\lambda=3$, the same result holds for
$n=13,17,27,30,45$. When $\lambda=4$, the same result holds for
$n=17,18,26,32,50,68,80$. Then there are infinitely many $k$ such
that $2,5,10,13,17,18,26,27,30,32,45,50,52,68,80$ are $k$-congruent
numbers.

In the following table, we give some integer solutions of
$(k^2-1)n=\alpha^4-\beta^4$ for $1<n\leq 10.$

\[\begin{tabular}{c|c}
\hline

 $n$ & $(k,x,y)$ \\

\hline

$2$ & $(11, 4, 2), (131, 14, 8), (181, 16, 2), (513, 34, 30), (573, 29, 15)$  \\

\hline

$3$ & $(9, 4, 2), (57, 10, 4), (521, 32, 22), (729, 37, 23)$ \\

\hline

$4$ & $(31, 8, 4), (59, 13, 11), (129, 18, 14), (161, 18, 6), (365, 31, 25),$ \\ &$(511, 32, 8), (545, 44, 40)$ \\

\hline

$5$ & $(2, 2, 1), (6, 4, 3), (7, 4, 2), (86, 16, 13), (390, 58, 57), (482, 38, 31),$ \\ &$ (487, 33, 3), (985, 47, 13)$\\

\hline

$6$ & $(69, 13, 1), (219, 34, 32), (319, 28, 8), (441, 37, 29)$ \\
\hline

$7$ &$(103, 22, 20), (519, 51, 47)$ \\
\hline

$8$ & $(33, 10, 6), (239, 26, 2), (481, 38, 22), (611, 73, 71), (781, 47, 1)$ \\
\hline

$9$ & $(649, 57, 51)$ \\
\hline

$10$ & $(3, 3, 1), (5, 4, 2), (35, 11, 7), (83, 17, 11), (365, 34, 8), (581, 76, 74),$ \\ &$ (773, 52, 34), (897, 54, 26)$ \\
\hline
\end{tabular}
\]

Hence, we raise the following problem

\begin{problem} For each $n\geq 2,$ whether there are infinitely many $k$ such that $n$ is a $k$-congruent number.
\end{problem}

By Tunnell's theorem \cite{Tunnell} about congruent numbers, we have

\begin{corollary} Assuming the validity of the BSD Conjecture for $E_m:~y^2=x^3-m^2x,$ the following statements are
equivalent:\\

$(1)$ For a given $k\geq 2$, $n$ is a $k$-congruent number.\\

$(2)$ If $(k^2-1)n$ is odd, then the number of triples of integers
$(x,y,z)$ satisfying $2x^2+y^2+8z^2=(k^2-1)n$ is equal to twice the
number of triples satisfying $2x^2+y^2+32z^2=(k^2-1)n$. If
$(k^2-1)n$ is even, then the number of triples of integers $(x,y,z)$
satisfying $8x^2+2y^2+16z^2=(k^2-1)n$ is equal to twice the number
of triples satisfying $8x^2+2y^2+64z^2=(k^2-1)n$.
\end{corollary}

\vskip 10pt
\section{Further Consideration on $E_{n,d}$}
In the proof of Theorem 1.6, we get a family of elliptic curves
\begin{equation}
E_{n,d}: y^2=x^3-\frac{3n^2+d^4}{3}x+\frac{(9n^2+2d^4)d^2}{27},
\end{equation}
which seems complicated, but there are some interesting things. Multiply 729 on both sides of (4.1), we have
\begin{equation}
E'_{n,d}: y^2=x^3-(81n^2+27d^4)x+27d^2(9n^2+2d^4).
\end{equation}
For $n\neq d^2,$ we have the following proposition.

\begin{proposition}If $n\neq d^2$, then $n$ is a $d$-congruent number.
\end{proposition}
\emph{Proof.} By some calculations, we find that the points
$Q=(-6d^2,27dn)$ and $R=(3d^2-9n,27dn)$ lie on $E'_{n,d}$. From the
group law, we get
\[[2]Q=\bigg(\frac{3(n^2+3d^4)(3n^2+d^4)}{4d^2n^2},\frac{-27(n^2+d^4)(d^8+4d^4n^2-n^4)}{8d^3n^3}\bigg).\]
Let the point $S$ be the intersection of the line, which goes through $Q$, $R$, and the elliptic curve $E'_{n,d}$, we have
\[S=\bigg(\frac{3(d^6-nd^4+7d^2n^2-3n^3)}{(n+d^2)^2},\frac{-27dn(-n+d^2)(d^4+3n^2)}{(n+d^2)^3}\bigg),\]
which leads to
 \[\begin{cases}\begin{split}
&a=\frac{2(d^4+n^2)d}{(n-d^2)(n+d^2)},\\
&b=\frac{(n-d^2)(n+d^2)}{2nd},\\
&c=\frac{n^4+6d^4n^2+d^8}{2(n-d^2)(n+d^2)dn},
\end{split}\end{cases}\]
it's easy to see that $a,b,c\in \mathbb{Q}^+$ when $n>d^2.$

From the point \[-S=[-1]S=\bigg(\frac{3(d^6-nd^4+7d^2n^2-3n^3)}{(n+d^2)^2},\frac{27dn(-n+d^2)(d^4+3n^2)}{(n+d^2)^3}\bigg),\]
we get
 \[\begin{cases}\begin{split}
&a=\frac{4dn^2}{(-n+d^2)(n+d^2)},\\
&b=\frac{n(-n+d^2)(n+d^2)}{(d^4+n^2)d},\\
&c=\frac{n(n^4+6d^4n^2+d^8)}{2(-n+d^2)(n+d^2)(d^4+n^2)d},
\end{split}\end{cases}\]
it's easy to see that $a,b,c\in \mathbb{Q}^+$ when $n<d^2.$

Therefore, for $n\neq d^2,$ all other positive integers are $d$-congruent numbers. \hfill $\Box$\\

For example, when $d=1,$ for $n\geq 2$, we have
 \[\begin{cases}\begin{split}
&a=\frac{2(n^2+1)}{(n-1)(n+1)},\\
&b=\frac{(n-1)(n+1)}{2n},\\
&c=\frac{n^4+6n^2+1}{2(n-1)(n+1)n},
\end{split}\end{cases}\]
i.e., all $n\geq 2$ are $1$-congruent numbers.

For a fixed $n$, from the proof of Proposition 4.1, we have

\begin{proposition}For any $n\in \mathbb{Z}^+$, except $d^2=n$, all other $d$ such that $n$ is a $d$-congruent number.
\end{proposition}

\begin{remark} The $k$-congruent number curve $E_{n,k}$ and the $d$-congruent number
curve $E_{n,d}$ have a big difference. The $j$-invariant of
$E_{n,k}$ is 1728, but the $j$-invariant of $E_{n,d}$ depends on
$n$.

\end{remark}

\vskip 10pt
\bibliographystyle{amsplain}

\end{document}